\documentclass[12pt]{article}
\usepackage{amssymb,amsmath,amsfonts}
\usepackage{fullpage}

\def\incl{\lhook\joinrel@>>>\!\!\!\!@>>>}
\newtheorem{thm}{Theorem}[section]

\begin{document}

\title{Regularity of monoids under Sch\"{u}tzenberger products}
\author{Firat Ate\c{s} \hspace{0.2cm} {\footnotesize and} \hspace{0.2cm} A.
Sinan \c{C}evik \\
{\footnotesize firat@balikesir.edu.tr} \hspace{0.2cm}{\footnotesize and} 
\hspace{0.2cm} {\footnotesize scevik@balikesir.edu.tr}}
\date{Department of Mathematics, Faculty of Science and Art, \\
Balikesir University, \c{C}agis Campus, 10145, Balikesir, Turkey \\
[0.3cm]}

\maketitle

\begin{abstract}
{\footnotesize In this paper we give a partial answer to the problem which is about the regularity of Sch\"{u}tzenberger products in semigroups asked by Gallagher in his thesis \cite[Problem 6.1.6]{Gallagher} and, also, we investigate the regularity for the new version of the Sch\"{u}tzenberger product which was defined in \cite{Ates}.} \\[0.5mm]
\footnotesize{{\it Keywords}: The Sch\"{u}tzenberger product, regularity.   } 
\end{abstract}

\normalsize

\section{Introduction and Preliminaries}

\hspace{0.5cm} 
We recall that a monoid $M$ is called \textit{regular} if, for every $a \in M$, there exists $b \in M$ such that $aba=a$ and $bab=b$ (or, equivalently, for the set of inverses of $a$ in $M$, that is, $a^{-1}=\{b\in B$ : $aba=a$ and $bab=b\}$, $M$ is regular if and only if, for all $a\in M$, the set $a^{-1}$ is not equal to the emptyset). In \cite[Problem 6.1.6]{Gallagher}, Gallagher asked whether there exists a classification for arbitrary semigroups $A$ and $B$ for which the Sch\"{u}tzenberger product $A\Diamond B$ is regular. In fact, before asking this problem, the question of the regularity of the wreath product of monoids was explained by Skornjakov (\cite{Skorn}). After that, in \cite{Nico}, it has been investigated the regular properties of semidirect and wreath products of
monoids. In this paper, to convience the above problem, we purpose to give a partial answer by defining necessary and sufficient conditions of the Sch\"{u}tzenberger product $A\Diamond B$ to be regular where both $A$ and $B$ are any monoids. Moreover, by giving a new version of the Sch\"{u}tzenberger product (\cite{Ates}), say $A\Diamond_v B$, we will present another result about this regularity problem.
\par
A generating and defining relation sets for the Sch\"{u}tzenberger product of arbitrary monoids have been defined in a joint paper written by Howie and Ruskuc (in \cite{Howie-Ruskuc}). Moreover, in \cite{Gallagher}, Gallagher defined the finitely generatability and finitely presentability of this product and then he left an open problem explained in the above paragraph.  
\par
Let $A$ and $B$ be any monoids with associated presentations $\wp _{A}=[X;R]$ and $\wp _{B}=[Y;S]$, respectively. Each paragraph at the rest of this section, we will recall definitions of some products which will be needed for the main results of this paper.
\par
Let $M=A\rtimes _{\theta }B$ be the corresponding semidirect products of these two monoids, where $\theta $ is a monoid homomorphism from $B$ to $End(A)$ such that, for every $a\in A$, $b_{1},b_{2}\in B$, $(a)\theta _{b_{1}b_{2}}=((a)\theta _{b_{2}})\theta _{b_{1}}$.
We recall that the elements of $M$ can be regarded as ordered pairs $(a,b)$, where $a\in A$, $b\in B$ with the multiplication given by 
$(a_{1},b_{1})(a_{2},b_{2})=(a_{1}(a_{2})\theta _{b_{1}},b_{1}b_{2}),$
and the monoids $A$ and $B$ are identified with the submonoids of $M$ having elements $(a,1_{B})$ and $(1_{A},b)$. For every $x\in X$ and $y\in Y,$ choose a word, denoted by $(x)\theta _{y}$, on $X$ such that $[(x)\theta _{y}]=[x]\theta _{\lbrack y]}$ as an element of $K$. To establish
notation, let us denote the relation $yx=(x)\theta _{y}y$ on $X\cup Y$ by $T_{yx}$ and write $T$ for the set of relations $T_{yx}$. Then, for any
choice of the words $(x)\theta _{y}$, $\wp _{M}=[X,Y \: ; \: R, S, T]$
is a standard monoid presentation for the semidirect product $M$.
\par
The cartesian product of $B$ copies of the monoid $A$ is denoted by $A^{\times B},$ while the corresponding direct
product is denoted by $A^{\oplus B}.$ One may think of $A^{\times B}$ as the set of all such functions from $B$ to $A$, and $A^{\oplus B}$ as the set all such functions $f$ having finite support, that is to say, having the property that $(x)f=1_{A}$ for all but finitely many $x$ in $B.$ The
unrestricted and restricted wreath products of the monoid $A$ by the monoid $B$, are the sets $A^{\times B}\times B$ and $A^{\oplus B}\times B$,
respectively, with the multiplication defined by $(f,b)(g,b^{^{\prime }})=(f\text{ }^{b}g,bb^{^{\prime }}),$
where $^{b}g:B\rightarrow A$ is defined by 
\begin{equation}  
\label{1}
(x)^{b}g=(xb)g,\quad (x\in B)
\end{equation}
such that $(xb)g$ has finite support. It is well known that both these wreath products are monoids with the identity $(\overline{1},1_{B}),$ where $x\overline{1}=1_{A}$ for all $x\in B$. (For more details on the definition and applications of restricted (unrestricted) wreath products, we can refer, for instance, \cite{Baumslag, Howie-Ruskuc, Meldrum, Ruskuc-1, Ruskuc-2}). We should note that, for having finite support, $B$ must be finite or groups.
\par
Now for a subset $P$ of $A\times B$ and $a\in A,$ $b\in B,$ we let define 
$$
Pb=\{(c,db)\; ; \; (c,d)\in P\}\text{ \ and \ }aP=\{(ac,d)\; ;\; (c,d)\in P\}.
$$
Then the Sch\"{u}tzenberger product of $A$ and $B$, denoted by $A\Diamond B$, is the set $A\times P(A\times B)\times B$ with the multiplication
$(a_{1},P_{1},b_{1})(a_{2},P_{2},b_{2})=(a_{1}a_{2},P_{1}b_{2}\cup
a_{1}P_{2},b_{1}b_{2}).$
Clearly $A\Diamond B$ is a monoid (\cite{Howie-Ruskuc}) with the identity $(1_{A}, \emptyset, 1_{B})$.

\section{Main Theorems}

\hspace{0.5cm} The following first theorem aims to give necessary and sufficient conditions for $A\Diamond B$ to be regular while both $A$ and $B$ are arbitrary monoids.
\begin{thm}
\label{theoremreg1}
Let $A$ and $B$ be any monoids. The product $A\Diamond B$ is regular if and only if
\begin{enumerate}
\item[$(i)$] $A$ and $B$ are regular,
\item[$(ii)$] for every $(a, P, b)\in A\Diamond B$, either 
$$
P=aP_{1}b=\bigcup_{(a_{1}, b_{1})\in P_{1}}\{(aa_{1}, b_{1} b)\} \quad or \quad P=caP_{1}bd=\bigcup_{(a_{1}, b_{1})\in P_{1}}\{(caa_{1}, b_{1}bd)\},
$$ 
where $P_{1}\subseteq A\times B$ and $c\in a^{-1}$, $d\in b^{-1}$.
\end{enumerate}
\end{thm}
\par
By (\ref{1}) and the definiton of Sch\"{u}tzenberger product, we can define a new version of the Sch\"{u}tzenberger product as follows. We note that the definition and some other properties of this product have been investigated in \cite{Ates}. 
\par
Let $A$ and $B$ be monoids. We recall that $A^{\oplus B}$ is
the set of all functions $f$ having finite support. For $P\subseteq
A^{\oplus B}\times B$ and $b\in B$, we define the set 
\begin{equation*}
Pb=\{(f,db);(f,d)\in P\}.
\end{equation*}
The new version of the Sch\"{u}tzenberger product of $A$ by $B$, denoted by $
A \Diamond _{v} B$, is the set $A^{\oplus B}\times P(A^{\oplus B}\times
B)\times B$ with the multiplication 
\begin{equation*}
(f,P_{1},b_{1})(g,P_{2},b_{2})=(f\text{ }^{b_{1}}g,P_{1}b_{2}\cup
P_{2},b_{1}b_{2}).
\end{equation*}
One can easily show that $A \Diamond _{v}B$ is a monoid with the identity $(
\overline{1},\emptyset ,1_{B})$, where $^{b_1}g$ is defined as in (\ref{1}). We should also note that, for having finite support, $B$ must be finite or groups.
\par
Thus another main result of this paper is the following.
\begin{thm}
\label{theoremreg2} 
Let $A$ be an arbitrary monoids and $B$ be a finite monoid or be a group. Then $A \Diamond _{v} B$ is regular if and only if 
\begin{enumerate}
\item[$(i)$] $A$ and $B$ are regular,
\item[$(ii)$] For every $x\in B$ and $f\in A^{\oplus B}$ there exist $e \in B$ such that $e^{2}=e$, with 
$$
(x)f\in A(xe)f.
$$
\item[$(iii)$] for every $(f, P, b)\in A \Diamond _{v} B$, either 
$$
P=P_{1}b=\bigcup_{(f_{1}, b_{1})\in P_{1}}\{(f_{1}, b_{1} b)\} \quad or \quad P=P_{1}bd=\bigcup_{(f_{1}, b_{1})\in P_{1}}\{(f_{1}, b_{1}bd)\},
$$ 
where $P_{1}\subseteq A^{\oplus B}\times B$ and $d\in b^{-1}$.
\end{enumerate}
\end{thm}

\section{Proofs}

{\bf Proof of Theorem \ref{theoremreg1}:}
Let us suppose that $A\Diamond B$ is regular. Thus, for $(a,\emptyset , b)\in $ $A\Diamond B$, there exists $(c,P,d)$ such that
\begin{eqnarray*}
(a,\emptyset , b) &=&(a, \emptyset , b)(c,P,d)(a, \emptyset , b)=(aca, aPb, bdb), \\
(c,P,d) &=&(c,P,d)(a, \emptyset ,b)(c,P,d)=(cac, Pbd \cup{caP}, dbd).
\end{eqnarray*}
Therefore we have $a=aca$, $c=cac$, $b=bdb$ and $d=dbd$. This implies that $(i)$ must hold. 
\par
By the assumption on the regularity of $A\Diamond B$, for $(a,P,b)\in $ $A\Diamond B$, we have $(c,P_{2},d)\in $ $A\Diamond B$ such that
$$
(a, P, b) = (a, P, b)(c, P_{2}, d)(a, P, b) \quad \mbox{and} \quad (c, P_{2}, d) = (c, P_{2}, d)(a, P, b)(c, P_{2}, d).
$$
Hence this gives us $a=aca$, $c=cac$, $b=bdb$, $d=dbd$, $P=Pdb\cup aP_{2}b \cup acP$ and $P_{2}=P_{2}bd \cup cPd\cup caP_{2}$. To show the second condition in theorem, let us suppose that $P\neq aP_{1}b$, for some $P_1 \subseteq A \times B$. Then there exists $(a_2,b_2)\in P$ such that $a_2\neq aa_{2}^{'}$ and $b_2\neq b_{2}^{'}b$ where $a_{2}^{'}\in A$ and $b_{2}^{'}\in B$.
Thus $P$ can not be equal to
$
Pdb\cup aP_{2}b\cup acP
$, for all $P_2 \subseteq A\times B$. This gives a contradiction with the regularity of $A\Diamond B$. In fact, when someone take $P=aP_{1}b$, the equalities 
\begin{eqnarray*}
Pdb\cup aP_{2}b\cup acP =  aP_{1}bdb\cup aP_{2}b\cup acaP_{1}b& = & aP_{1}b\cup aP_{2}b\cup aP_{1}b \\
 & = & aP_1b \quad \mbox{by choosing } P_2= caP_{1}bd \\
 & = & P
\end{eqnarray*}
and
\begin{eqnarray*}
P_{2}bd \cup cPd\cup caP_{2} & = & P_{2}bd \cup caP_1bd \cup caP_{2} \\
 & = & caP_{1}bdbd \cup caP_1bd \cup cacaP_{1}bd \\
 & & \mbox{by choosing } P_2= caP_{1}bd \\
 & = & caP_{1}bd \cup caP_1bd \cup caP_{1}bd = caP_1bd = P_2
\end{eqnarray*}
hold. We note that, by applying similar discussions as above for the case $P=caP_{1}bd$ in theorem, where $P_{1}\subseteq A\times B$ and $c\in a^{-1}$, it is seen that condition $(ii)$ must hold.
\par
For the converse part of the proof, let $(a, P, b)\in A\Diamond B$. Thus we definitely have $c\in A$ and $d\in B$ such that $c\in a^{-1}$ and $d\in b^{-1}$. Now let us consider the union of sets
$$
Pdb \cup aP_2b \cup acP \quad \mbox{and} \quad P_2bd \cup cPd \cup caP_2.
$$
At this stage, by $P=aP_{1}b$, if we choose $P_{2}=caP_{1}bd \subseteq A\times B$, then we get  
$$
Pdb\cup aP_{2}b\cup acP = aP_{1}b=P \quad \mbox{and} \quad P_{2}bd\cup cPd\cup caP_{2} = caP_{1}bd= P_{2}.
$$
As a result of this, for every $(a, P, b)\in A\Diamond B$, there exists $(c, P_2, d)\in A\Diamond B$ such that 
\begin{eqnarray*}
(a, P, b)(c, P_{2}, d)(a, P, b)=(aca, Pdb\cup aP_{2}b\cup acP,bdb)=(a, P, b), \\
(c, P_{2}, d)(a, P, b)(c, P_{2}, d) = (cac,P_{2}bd\cup cPd\cup caP_{2}, dbd) = (c, P_{2}, d).
\end{eqnarray*}
In addition, by applying similar above arguments for the case $P=caP_{1}bd$ in theorem, where $P_{1}\subseteq A\times B$ and $c\in a^{-1}$, the proof of the regularity of $A \Diamond B$ is completed.
\par
Hence the result. {$\square$}\\[0.5cm]
{\bf Proof of Theorem \ref{theoremreg2}:}
{\normalsize Let us suppose that $A \Diamond _{v} B$ is regular. Thus, for $(f,(1_{A},1_{B}),b)\in A \Diamond _{v} B$, there exists $(g
,P,d)\in A \Diamond _{v} B$ such that 
\begin{eqnarray*}
(f,(1_{A},1_{B}),b) &=&(f,(1_{A},1_{B}),b)
(g,P,d)(f,(1_{A},1_{B}),b), \\
(g,P,d) &=&(g,P,d)(f,(1_{A},1_{B}),b)(g,P,d).
\end{eqnarray*}
We then have $b=bdb$ and $d=dbd$. If we choose $b=1$ then we have $bd=1$. Therefore we have $f=fgf$ and $g=gfg$. This implies that both $B$ and $A^{\oplus B}$ are regular. Since $A^{\oplus B}$ denotes the direct product of $B$ copies of $A$, it is easy to see that if $A^{\oplus B}$ is regular, then $A$ is regular. This gives condition $(i)$. 
\par
By the assumption, for every $(f,P,b)\in A \Diamond _{v} B$, we have $(g,P_{2},d)\in A \Diamond _{v} B$ such that 
\begin{eqnarray*}
(f,P,b) &=&(f,P,b)(g,P_{2},d)(f,P,b)=(f\text{ }^{b}g\text{ }^{bd}f,Pdb \cup P_{2}b \cup P,bdb), \\
(g,P_{2},d) &=&(g,P_{2},d)(f,P,b)(g,P_{2},d)=(g\text{ }^{d}f\text{ }^{db}g,P_{2}bd \cup Pd \cup P_{2} ,dbd). 
\end{eqnarray*}
Hence, by equating the components, we get $f=f\: ^{b}g \: ^{bd}f$, $g=g\: ^{d}f\: ^{db}g$, $b=bdb$, $d=dbd$, $P=Pdb\cup P_{2}b\cup P$ and $P_{2}=P_{2}bd\cup Pd\cup P_{2}$. These show that, for every $x\in B$,
\begin{eqnarray*} 
(x)f &=&(x)f\text{ }(x)^{b}g\text{ }(x)^{bd}f=(x)f\text{ }(xb)g\text{ }(xbd)f\in A(xbd)f.
\end{eqnarray*}
If we take $e=bd$, then condition $(ii)$ becomes true. In addition, by using the facts $b=bdb$, $d=dbd$, $P=Pdb\cup P_{2}b\cup P$ and $P_{2}=P_{2}bd\cup Pd\cup P_{2}$, for every $(f, P, b)\in A \Diamond _{v} B$, and by applying similar arguments given in the proof of Theorem \ref{theoremreg1}, we get 
$$
\mbox{either} \quad P=P_{1}b \quad \mbox{or} \quad P=P_{1}bd,
$$ 
where $P_{1} \subseteq A^{\oplus B}\times B$ and $d\in b^{-1}$. Therefore condition $(iii)$ must hold. 
\par
Conversely, let us suppose that the monoids $A$ and $B$ satisfy conditions $(i)$, $(ii)$ and $(iii)$. For $x,b,d\in B$ and $f,g\in A^{\oplus B},$
we let consider 
$$
(x)f\text{ }(x)^{b}g\text{ }(x)^{bd}f,
$$
where $dbd=d.$ By condition $(ii)$, for $a\in A$, we have $(x)f=a(xbd)f$ where $bd=e$. Thus 
\begin{equation}  \label{11}
(x)f\text{ }(x)^{b}g\text{ }(x)^{bd}f=a(xbd)f\text{ }(x)^{b}g\text{ } (x)^{bd}f=a(x)^{bd}f\text{ }(x)^{b}g\text{ }(x)^{bd}f.
\end{equation}
Since $A$ is regular, $A^{\oplus B}$ is regular \cite{Nico}. Thus we can choose $g= \: ^{d}v$ $(v \in A^{\oplus B})$ such that $fvf=f$ and $vfv=v.$ Hence the last term in (\ref{11}) will be equal to 
$$
a(x)^{bd}f\text{ }(x)^{bd}v\text{ }(x)^{bd}f=a(x)^{bd}(fvf)=a(x)^{bd}f=(x)f.
$$
This implies that $f=f \: ^{b}g \: ^{bd}f$. On the other hand, by similar procedure as above, we obtain 
$$
g\text{ }^{d}f\text{ }^{db}g=\text{ }^{d}v\text{ }^{d}f\text{ }^{dbd}v=\text{
}^{d}v\text{ }^{d}f\text{ }^{d}v=\text{ }^{d}(vfv)=\text{ }^{d}v=g.
$$
Moreover, by condition $(iii)$, we have $P=P_{1}b$ or $P = P_1bd$, where $P_{1}\subseteq A^{\oplus B}\times B$. For the next stage of proof, we will only consider $P = P_1b$ since similar progress can be applied for the other value of $P$. Therefore there exists a subset $P_{2}=P_{1}bd$ of $A^{\oplus B}\times B$ such that 
\begin{eqnarray*}
Pdb\cup P_{2}b\cup P = P_{1}bdb\cup P_{1}bdb\cup P_{1}b=P_{1}b\cup
P_{1}b\cup P_{1}b=P_{1}b=P, \\
P_{2}bd\cup Pd\cup P_{2} = P_{1}bdbd\cup P_{1}bd\cup P_{1}bd=P_{1}bd\cup
P_{1}bd\cup P_{1}bd=P_{1}bd=P_{2}.
\end{eqnarray*}
As a result of these above procedure, for every $(f, P, b)\in A\Diamond _{v}B$, there exists $(g, P_{2}, d)\in A\Diamond _{v}B$ such that 
\begin{eqnarray*}
(f, P, b)(g, P_{2}, d)(f, P, b)=(f \: ^{b}g \: ^{bd}f, Pdb\cup P_{2}b\cup P, bdb)=(f, P, b), \\
(g, P_{2}, d)(f, P, b)(g, P_{2}, d) = (g \: ^{d}f \: ^{db}g,P_{2}bd\cup Pd\cup P_{2}, dbd) = (g, P_{2}, d).
\end{eqnarray*}
Hence the result. {$\square$\par\bigskip}

\vskip 1cm 

\begin{center}
\textit{The corresponding addresses for the authors:} \\[0pt]
{\footnotesize {Firat Ate\c{s} \hspace{0.3cm} and \hspace{0.3cm} A. Sinan 
\c{C}evik \\[0pt]
Balikesir Universitesi,\\[0pt]
Fen-Edebiyat Fakultesi, \\[0pt]
Matematik Bolumu, Cagis Kampusu, \\[0pt]
10145 Balikesir/TURKEY }\\[0pt]
e-mails: \textit{firat@balikesir.edu.tr} \hspace{0.3cm} and \hspace{0.3cm} 
\textit{scevik@balikesir.edu.tr} }
\end{center}
\end{document}